# Constructing Even Order Magic Squares By Consecutive Numbering

A. M. IBRAHIM and H. M. JIBRIL
Ahmadu Bello University
Zaria, Nigeria.

ABDULLAHI UMAR
Sultan Qaboos University
Al-Khod, P. C. 123, Muscat, OMAN

**Abstract**

The aim of this note is to introduce fastest new general methods for the construction of double and single even order magic squares. As in [5], the method for double even order magic squares is fairly straight-forward but some adjustments are necessary for the single even order magic squares.

## 1. Introduction

For centuries magic squares have fascinated specialists and non-specialists alike, and still continue to do so. The study of magic squares is basically concerned with two problems: finding general methods for constructing all $n$ by $n$ magic squares, or constructing all magic squares of a given class, and enumerating how many distinct (excluding reflections and rotations) magic squares there are for a given order $n$. The latter seems to be the more difficult problem as very little progress has been made. In fact, there is a unique 3 by 3 magic square, exactly 880 (eight-hundred and eighty) 4 by 4 magic squares. In [6], it is claimed that the number of 5 by 5 magic squares was computed by R. Schroeppel in 1973 as 275305224 while Pinn and Wieczerkowski estimated the number of 6 by 6 magic squares to be $(1.7745 \pm 0.0016) \times 10^{19}$ using Monte Carlo simulation and methods from statistical mechanics.

There are now several methods for constructing magic squares of any order (see references [1] – [5] for example). The fastest and general method for constructing odd order magic squares by consecutive numbering was found by the Middle

Easterners. This method appears fully developed in Persian manuscripts of the 13[th] century [3]. In this note we describe constructions of even order magic squares by consecutive numbering, which we believe are new. This construction is a slight modification of the method described in [5] but is far better. However, as observed in [5], it is necessary to classify even order magic squares into doubly and singly-even and note that the squares of the latter type are more difficult to construct.

## 2. Basic facts and definitions

A *primitive magic square* (referred to as a magic square in what follows) of *order n* is a square array consisting of the $n^2$ distinct numbers *1, 2, 3, … , $n^2$* such that the sum of the entries in each row, column and main diagonals add up to the same total, $n(n^2 + 1)/2$.

A *doubly-even order* magic square is one whose order is divisible by 4.

A *singly-even order* magic square is one whose order is divisible by 2 but not divisible by 4.

A pair of natural numbers *(a, b)* is *complementary* if $a + b = n^2 + 1$, where *n* is the order of the magic square and, *a* and *b* are called *complements* of each other. There are $n^2/2$ distinct complementary pairs in an even order magic square.

A magic square is *parallel* if all its complementary pairs are placed parallel to one another.

An even order magic square is *associated* if all its complementary pairs are equidistant from the centre of symmetry of the square, i. e., the point where all the lines of symmetry meet.

A magic square is *mixed* if it is neither parallel nor associated.

## 3. The construction of doubly-even order magic squares

Since a magic square is a sort of symmetrical or balanced square we may therefore impose a simple symmetry on the even-odd number distribution. First, let's consider the standard 4-quadrant square given by Fig I(a), where Quadrant IV is a vertical translation of Quadrant I and the left half of the square is a reflection (along the vertical line of symmetry) of the right half. On each 2 by 2 subsquare within any quadrant we impose a simple even-odd symmetry given by Figure I(b) or its rotation (through $90^0$ or $270^0$).

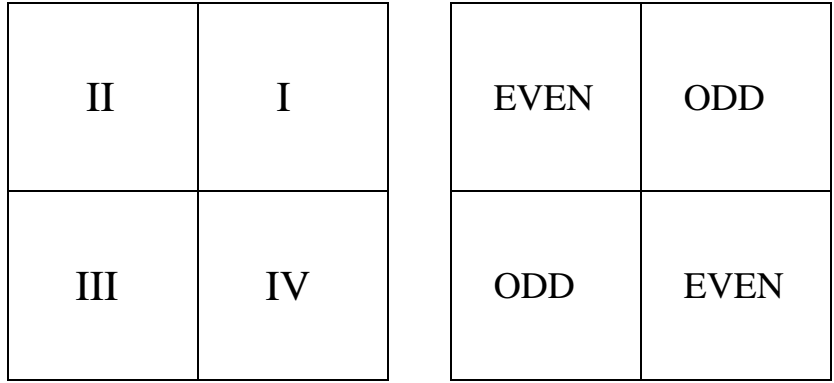

Fig. I (a) Standard 4-quadrant square, (b) Even-odd symmetric *2 by 2* subsquare

We are going to describe the detailed construction in three steps.

STEP 1. Let $p = n^2/2$ and $m = n/2$. For $k = 1, 2, ..., m$ consider the $n$ complementary pairs:

$((k-1)n+1, 2p-(k-1)n), ((k-1)n+2, 2p-(k-1)n-1), ..., (kn-1, 2p-kn+2),$
$(kn, 2p-kn+1);$

and rearrange them into new (noncomplementary) pairs as follows:

$((k-1)n+1, 2p-kn+1), ((k-1)n+2, 2p-kn+2), ..., (kn-1, 2p-(k-1)n-1),$
$(kn, 2p-(k-1)n).$

STEP 2. Place these pairs consecutively (with respect to the first numbers in each pair) in the *k*th and (*n* – *k* + 1)th-columns from top to bottom for odd *k*; and from bottom to top for even *k*. Notice that the resulting square is associated (though not yet magic), and so each main diagonal consists of *m* complementary pairs. The sum of the *r*th-row is

$$\sum_{k-odd}[((k-1)n+r+(2p-kn+r)] + \sum_{k-even}[(kn-(r-1))+(2p-(k-1)n-(r-1)]$$
$$= (2p-n+2r)(n/4) + (2p+n-2r+2)(n/4) = (4p+2)(n/4) = m(2p+1).$$

STEP 3. For rows— *r = 2, 4, ... , m, m + 1, m + 3, ... , n – 1* we swap (row-wise) the pairs in each symmetric pair of columns to achieve the even-odd symmetry given by Fig I(b). Note that this has no effect on the row sums or the associatedness of the square, however, we now have the extra advantage of fixing the sum of each column to the required sum of $n(n^2 + 1)/2 = m(2p + 1)$. In fact, the kth-column now consists of four subsequences:

$$\sum_{i=1}^{n/4}[(k-1)n+(2i-1)] + \sum_{i=1}^{n/4}[2p-kn+2i)] + \sum_{i=1}^{n/4}[2p-(k-1)n-(2i-1)] + \sum_{i=1}^{n/4}[kn-(2i-2)]$$

$$= (4p+2)(n/4)$$

$$= m(2p+1).$$

The construction of an *8 by 8* square is shown in Figures II (a) and (b), while the general result for an *n by n* associated magic square is shown in Figure III.

| 1 | 16 | 17 | 32 | 40 | 41 | 56 | 57 |
|---|----|----|----|----|----|----|----|
| 2 | 15 | 18 | 31 | 39 | 42 | 55 | 58 |
| 3 | 14 | 19 | 30 | 38 | 43 | 54 | 59 |
| 4 | 13 | 20 | 29 | 37 | 44 | 53 | 60 |
| 5 | 12 | 21 | 28 | 36 | 45 | 52 | 61 |
| 6 | 11 | 22 | 27 | 35 | 46 | 51 | 62 |
| 7 | 10 | 23 | 26 | 34 | 47 | 50 | 63 |
| 8 | 9  | 24 | 25 | 33 | 48 | 49 | 64 |

| 1  | 16 | 17 | 32 | 40 | 41 | 56 | 57 |
|----|----|----|----|----|----|----|----|
| 58 | 55 | 42 | 39 | 31 | 18 | 15 | 2  |
| 3  | 14 | 19 | 30 | 38 | 43 | 54 | 59 |
| 60 | 53 | 44 | 37 | 29 | 20 | 13 | 4  |
| 61 | 52 | 45 | 36 | 28 | 21 | 12 | 5  |
| 6  | 11 | 22 | 27 | 35 | 46 | 51 | 62 |
| 63 | 50 | 47 | 34 | 26 | 23 | 10 | 7  |
| 8  | 9  | 24 | 25 | 33 | 48 | 49 | 64 |

Fig. II (a) Constructing an *8 by 8* square,       (b) A complete *8 by 8* square

**The construction of single even order magic square.**

As remarked in the introduction, these are the most difficult magic squares to construct, and general methods are rare. The process is in many ways similar to that of the double even case described above. Ignoring the first and last rows we fill the $(n-2)\,by\,n$ subsquare as in the double even case described above, using only the first *(p – n)/2* and the last *(p – n)/2* complementary pairs, i. e., omitting the *n* middle complementary pairs. The detailed description is as follows:

STEP 1. Let $p = n^2/2$ and $m = n/2$. For $k = 1,2,...,m-1$ consider the *n – 2* complementary pairs:

$((k-1)(n-2)+1, 2p-(k-1)(n-2), ((k-1)(n-2)+2, 2p-(k-1)(n-2)n-1),...,$

$(k(n-2)-1, 2p-k(n-2)+2), (k(n-2), 2p-k(n-2)+1);$

and rearrange them into new (noncomplementary) pairs as follows:

$((k-1)(n-2)+1, 2p-k(n-2)+1), ((k-1)(n-2)+2, 2p-k(n-2)+2),...,$

$(k(n-2)-1, 2p-(k-1)(n-2)-1), (k(n-2), 2p-(k-1)(n-2)).$

For $k = m$ the complementary pairs are retained, i. e., not to be rearranged.

STEP 2. Place these pairs consecutively (with respect to the smaller numbers in each pair) in the $k$th and $(n – k + 1)$th-columns from top to bottom for odd $k$; and from bottom to top for even $k$.

STEP 3. For rows— $r = 2, 4, ... , m, m + 1, m + 3, ... , n – 1$ we swap (row-wise) the pairs in each symmetric pair of columns. Note that this has no effect on the row sums, however, we now have the extra advantage of fixing the sum of each column to $n(n^2 + 1)/2 – (n^2 + 1) = (m – 1)(2p + 1)$.

STEP 4. Before we describe how to fill the outermost rows (R1 and R$n$), we first define a finite subsequence of the natural numbers. Let $\{a_j\}$ be a subsequence, with common difference $d = 1$ and whose first and last terms are given by

$$a_j : p – n + 1, ... , p + n.$$

In fact, $\{a_j\}$ consists of the middle $n$ complementary pairs omitted thus far, and has $2n$ terms.

The two outermost rows R1 and R$n$ are filled as shown in Figures IVa and IVb. There are $n$ cells to be filled from the left numbered 1 to $n$. Complementary pairs are placed vertically.

It is routine to verify that this construction yields a mixed magic square. We exemplify the process by constructing a *10 by 10* magic square in Figures V (a), (b) and (c).

**Concluding remarks.**

The main advantage of these methods of constructing even order magic squares over the numerous known methods is that, they can be regarded as 'continuous' or 'consecutive'. Thus, making them the fastest methods by analogy with the fastest method for constructing odd order magic squares. For the doubly-even case the description of the method is as follows:

'Beginning with the top left corner cell, place consecutively, the numbers $1, 2, \ldots, n$ from the first to the last row (one number in each row) alternating from left to right to left … in the pair of outermost columns. However, $n/2$ and $(n/2)+1$ are to be placed in the same column, after which you resume the alternation from right to left to right …. Next place the numbers $n+1, n+2, \ldots, 2n$ from the last row to the first row in

| 1 | 2n | 2n+1 | 4n | … | p | p+n | … | 2p-3n | 2p-3n+1 | 2p-n | 2p-n+1 |
|---|---|---|---|---|---|---|---|---|---|---|---|
| 2p-n+2 | 2p-n-1 | 2p-3n+2 | 2p-3n-1 | … | p+n-1 | p-1 | … | 4n-1 | 2n+2 | 2n-1 | 2 |
| 3 | 2n-2 | 2n+3 | 4n-2 | … | p-2 | p+n-2 | … | 2p-3n-2 | 2p-3n+3 | 2p-n-2 | 2p-n+3 |
| 2p-n+4 | 2p-n-3 | 2p-3n+4 | 2p-3n-3 | … | p+n-3 | p-3 | … | 4n-3 | 2n+4 | 2n-3 | 4 |
| . | . | . | . |  | . | . |  | . | . | . | . |
| . | . | . | . |  | . | . |  | . | . | . | . |
| . | . | . | . |  | . | . |  | . | . | . | . |
| m-1 | 3m+2 | 5m-1 | 7m+2 | … | p-m+2 | p+m+2 | … | 2p-7m+2 | 2p-5m-1 | 2p-3m+2 | 2p-m-1 |
| 2p-m | 2p-3m+1 | 2p-5m | 2p-7m+1 | … | p+m+1 | p-m+1 | … | 7m+1 | 5m | 3m+1 | m |
| 2p-m+1 | 2p-3m | 2p-5m+1 | 2p-7m | … | p+m | p-m | … | 7m | 5m+1 | 3m | m+1 |
| m+2 | 3m-1 | 5m+2 | 7m-1 | … | p-m-1 | p+m-1 | … | 2p-7m-1 | 2p-5m+2 | 2p-3m-1 | 2p-m+2 |
| . | . | . | . |  | . | . |  | . | . | . | . |
| . | . | . | . |  | . | . |  | . | . | . | . |
| . | . | . | . |  | . | . |  | . | . | . | . |
| 2p-3 | 2p-2n+4 | 2p-2n-3 | 2p-4n+4 | … | p+4 | p-n+4 | … | 3n+4 | 3n-3 | n+4 | n-3 |
| n-2 | n+3 | 3n-2 | 3n+3 | … | p-n+3 | p+3 | … | 2p-4n+3 | 2p-2n-2 | 2p-2n+3 | 2p-2 |
| 2p-1 | 2p-2n+2 | 2p-2n-1 | 2p-4n+2 | … | p+2 | p-n+2 | … | 3n+2 | 3n-1 | n+2 | n-1 |
| n | n+1 | 3n | 3n+1 | … | p-n+1 | p+1 | … | 2p-4n+1 | 2p-2n | 2p-2n+1 | 2p |

Fig. III A 'complete' *n* by *n* double even order associated magic square

| Position | 1 | 2 | 3 | 4 | 5 | ... | m+1 | m+2 | m+3 | m+4 | m+5 | m+6 | m+7 | ... | n |
|---|---|---|---|---|---|---|---|---|---|---|---|---|---|---|---|
| Number($R_1$) | | $a_1$ | | $a_3$ | | ... | $a_m$ | | | $a_{m+3}$ | | $a_{m+5}$ | | | |
| Number($R_n$) | | | $a_2$ | | $a_4$ | ... | | $a_{m+1}$ | $a_{m+2}$ | | $a_{m+4}$ | | $a_{m+6}$ | ... | $a_{n-1}$ |

Fig. IV (a) Filling the outermost rows of an *n by n* single even order magic square.

| Position | 1 | 2 | 3 | 4 | 5 | ... | m+1 | m+2 | m+3 | m+4 | m+5 | m+6 | m+7 | ... | n |
|---|---|---|---|---|---|---|---|---|---|---|---|---|---|---|---|
| Number($R_1$) | $a_{n+1}$ | $a_1$ | | $a_3$ | | ... | $a_m$ | | | $a_{m+3}$ | | $a_{m+5}$ | | | $a_{n+2}$ |
| Number($R_n$) | $a_n$ | | $a_2$ | | $a_4$ | ... | | $a_{m+1}$ | $a_{m+2}$ | | $a_{m+4}$ | | $a_{m+6}$ | ... | $a_{n-1}$ |

Fig. IV (b) The complete outermost rows of an *n by n* single even order magic square.

The blank cells consist of the complements of the numbers that appear above or below as the case may be.

|   |    |    |    |    |    |    |    |    |     |
|---|----|----|----|----|----|----|----|----|-----|
| 1 | 16 | 17 | 32 | 33 | 68 | 76 | 77 | 92 | 93  |
| 2 | 15 | 18 | 31 | 34 | 67 | 75 | 78 | 91 | 94  |
| 3 | 14 | 19 | 30 | 35 | 66 | 74 | 79 | 90 | 95  |
| 4 | 13 | 20 | 29 | 36 | 65 | 73 | 80 | 89 | 96  |
| 5 | 12 | 21 | 28 | 37 | 64 | 72 | 81 | 88 | 97  |
| 6 | 11 | 22 | 27 | 38 | 63 | 71 | 82 | 87 | 98  |
| 7 | 10 | 23 | 26 | 39 | 62 | 70 | 83 | 86 | 99  |
| 8 | 9  | 24 | 25 | 40 | 61 | 69 | 84 | 85 | 100 |
|   |    |    |    |    |    |    |    |    |     |

Fig. V (a) Constructing a *10 by 10* magic square with the first and last rows omitted

|    |    |    |    |    |    |    |    |    |     |
|----|----|----|----|----|----|----|----|----|-----|
| 1  | 16 | 17 | 32 | 33 | 68 | 76 | 77 | 92 | 93  |
| 94 | 91 | 78 | 75 | 67 | 34 | 31 | 18 | 15 | 2   |
| 3  | 14 | 19 | 30 | 35 | 66 | 74 | 79 | 90 | 95  |
| 96 | 89 | 80 | 73 | 65 | 36 | 29 | 20 | 13 | 4   |
| 97 | 88 | 81 | 72 | 64 | 37 | 28 | 21 | 12 | 5   |
| 6  | 11 | 22 | 27 | 38 | 63 | 71 | 82 | 87 | 98  |
| 99 | 86 | 83 | 70 | 62 | 39 | 26 | 23 | 10 | 7   |
| 8  | 9  | 24 | 25 | 40 | 61 | 69 | 84 | 85 | 100 |
|    |    |    |    |    |    |    |    |    |     |

Fig V (b)

| 51 | 41 | 59 | 43 | 57 | 45 | 55 | 54 | 48 | 52  |
|----|----|----|----|----|----|----|----|----|-----|
| 1  | 16 | 17 | 32 | 33 | 68 | 76 | 77 | 92 | 93  |
| 94 | 91 | 78 | 75 | 67 | 34 | 31 | 18 | 15 | 2   |
| 3  | 14 | 19 | 30 | 35 | 66 | 74 | 79 | 90 | 95  |
| 96 | 89 | 80 | 73 | 65 | 36 | 29 | 20 | 13 | 4   |
| 97 | 88 | 81 | 72 | 64 | 37 | 28 | 21 | 12 | 5   |
| 6  | 11 | 22 | 27 | 38 | 63 | 71 | 82 | 87 | 98  |
| 99 | 86 | 83 | 70 | 62 | 39 | 26 | 23 | 10 | 7   |
| 8  | 9  | 24 | 25 | 40 | 61 | 69 | 84 | 85 | 100 |
| 50 | 60 | 42 | 58 | 44 | 56 | 46 | 47 | 53 | 49  |

Fig. V (c) A complete *10 by 10* magic square exactly the same fashion as above in the pair of next outermost columns, noting that $3n/2$ and $(3n/2)+1$ are (vertically) adjacent, i. e., in the same column. The process continues until the pair of innermost columns is reached and $p$ is placed in the top cell of the left hand column of the innermost pair. Now place $p+1$ in the bottom cell of the right hand column of the innermost pair, and continue consecutive placement but reversing the whole process. That is, place the numbers: $p+1, p+2,...,p+n;$ $p+n+1, p+n+2,...,p+2n;$ ... ; $2p-n+1, 2p-n+2,...2p$ alternately from right (left) to left (right) to right (left) …, from the last (first) to the first (last) and moving outwards from the pair of innermost columns to the pair of outermost columns. Note that the process ends at the bottom right corner cell.'

And for the singly-even case the description of the method is as follows:

'Beginning with second cell from the top left corner, place consecutively, the numbers $1,2,...,n-2$ from the second to the last but one row (one number in each row) alternating from left to right to left … in the pair of outermost columns. However, $(n-2)/2$ and $(n-2)/2+1$ are to be placed in the same column, after which you resume the alternation from right to left to right …. Next place the numbers $n-1, n,...,2(n-2)$ from the last but one row to the second row in exactly the same fashion as above in the pair of next outermost columns, noting that $3(n-2)/2$ and $3(n-2)/2+1$ are adjacent. The process continues until the pair of innermost columns is reached and $q=(n^2-2n)/2$ is placed in the last but one cell of the left hand column of the innermost pair. Now we move into the pair of outermost (first and last) rows. Beginning with second cell (from the left) of the first row place consecutively, the numbers $q+1, q+2,...,q+n-1$ alternately from top to bottom (one in each colum, except the first) moving from the second to the last column. The top-bottom alternation pauses momentarily between $q+m+1$ and $q+m+2$, i. e., $q+m+1$ and $q+m+2$ are placed (horizontally) adjacent. The top-bottom alternation resumes immediately, until $q+n-1$ is placed in the bottom right hand corner. Then place $q+n$ in the bottom left hand corner, $q+n+1$ in the top left hand corner and $q+n+2$ in the top right hand corner. Next, place consecutively, the numbers $q+n+3, q+n+4,...,q+2n$ in the blank cells of the first and last row, alternately from the last but one column to the second column, until $q+2n$ is placed in the bottom cell of the second column. Now that the first and last rows are filled we go back into the square and place $q+2n+1$ adjacent (to the right) of $q$ and continue to fill the blank cells consecutively, from the innermost pair of columns to the outermost pair of columns, by the reverse of the original process. However, there is a break between $q+3n-2$ and $q+3n-1$, with $q+3n-1$ placed in the bottom (and not the top) blank cell of the next column and the process continues to completion.'

E-mail: aumarh@squ.edu.om